\documentclass[11pt]{article}
\usepackage{amsthm,amssymb,amscd}
\theoremstyle{plain}

\newtheorem{thm}{Theorem}[section]
\newcommand{\bb}[1]{\mbox{$\mathbb{#1}$}}

\title{On the comparison of different notions\\ of geometric categories}
\author{J\"{o}rg Sch\"{u}rmann\thanks{Westf. Wilhelms-Universit\"{a}t, SFB 478 
"Geometrische Strukturen in der Mathematik", Hittorfstr.27, 48149 M\"{u}nster,
Germany,
E-mail: jschuerm@math.uni-muenster.de}}
\date{ }

\begin{document}
\bibliographystyle{plain}

\maketitle

\begin{abstract}
In this paper we explain our notion of a 'Nash geometric 
category', which allows an easy comparison between the following different
axiomatic notions of geometric categories:
\begin{enumerate}
\item The o-minimal structures $\mathfrak{S}$ on the real field 
$(\bb{R},+,\cdot)$,
as defined by van den Dries \cite{vDr}.
\item The analytic geometric categories of van den Dries and Miller 
\cite{vDrMi}.
\item The $\mathfrak{X}$-sets of Shiota \cite{Shiota}.
\end{enumerate}
\end{abstract}

\section*{Introduction}  
In the last years different axiomatic generalizations of the theory of 
semialge\-braic- and subanalytic sets where developed.\\

First we have the theory of o-minimal structures, and especially of   
o-minimal structures on the real field 
$(\bb{R},+,\cdot)$, as defined by van den Dries and explained in his
book \cite{vDr}. This theory is an abstraction and generalization
of the theory of semialgebraic sets of affine spaces.\\

Another abstraction and generalization
of the theory of subanalytic sets of real analytic manifolds
is the theory of analytic geometric categories as introduced and developed
by van den Dries and Miller \cite{vDrMi}. This paper was motivated by the
work and the applications of Schmid and Vilonen about characteristic cycles 
of constructible complexes of sheaves \cite{SchVi} (as remarked by 
van den Dries and Miller
in the beginning of the introduction to \cite{vDrMi}).
Moreover, van den Dries and Miller explain a one to one correspondence between 
their analytic geometric categories and those 
o-minimal structures on the real field 
$(\bb{R},+,\cdot)$, which contain all globally subanalytic subsets.
They use this correspondence to transfer results from the o-minimal context
into the context of analytic geometric categories.\\

Finally, Shiota introduced another axiomatic notion
of $\mathfrak{X}$-sets of real affine spaces in his book \cite{Shiota}.
This theory is a simultaneously generalization of the theory of
semialgebraic- and subanalytic subsets of real affine spaces.\\

There are some results (e.g. results about curve selection, dimension,
Whitney stratifications and triangulability), which are quite similarly
in these categories (compare also with the article of Teissier \cite{Te}).
But there are also some important specific results, which are not yet
explained in all these different axiomatic generalizations:

\begin{enumerate}
\item The notion of a polynomially bounded o-minimal structure on the real 
field, related to the Lojasiewicz inequality, and the corresponding
growth dichotomy result (compare \cite[p.510/511]{vDrMi}).
\item The generic triviality results for o-minimal structures 
(\cite[4.11]{vDrMi} and \cite[chap. 9]{vDr}).
\item The results of van den Dries (\cite[chap. 10]{vDr}) about definable
spaces and quotients.  
\item The description of closed definable sets as the zero-set of
a definable $C^{p}$-function ($1\leq p< \infty$), as given by 
\cite[1.20,d.19]{vDrMi}.
\item The results in \cite{Shiota} about:\\
- Triangulability of definable functions and maps.\\
- Uniqueness results about this type of triangulations.\\
- $\mathfrak{X}$-versions of Thom's first and second isotopy lemmas.
\end{enumerate}

We introduced in our paper \cite{Sch} a variant of the 
analytic geometric categories of van den Dries and Miller \cite{vDrMi}, 
by restricting ourselves to 
(analytic) Nash manifolds. In this way, every
o-minimal structure $\mathfrak{S}$ on the real field $(\bb{R},+,\cdot)$ 
corresponds uniquely to a category, which we call a Nash geometric category.
This notion allows an easy comparison (contrary to a remark in 
\cite[p.498]{vDrMi}) between the different notions of 'geometric 
categories' as before. Moreover, it can easily be used to extend
the above results to all these different 'geometric categories'.\\

But we used this result in \cite{Sch} for the development
of the theory of constructible sheaves in these 'geometric categories'.
We recall in this paper our comparison result and its proof,
since it is not related to this abstract sheaf theory
(so this abstract language is not used in this paper).
Some people explained to me, that this should be useful.
Let us now recall the basic notions.\\

An o-minimal structure $\mathfrak{S}$ on $\bb{R}$ is a 
sequence $\mathfrak{S}_{n}$ ($n\in \bb{N}$) such that for each n:
\begin{enumerate}
\item $\mathfrak{S}_{n}$ is a boolean algebra of subsets of $\bb{R}^{n}$.
\item $A\in \mathfrak{S}_{n} \Rightarrow A\times\bb{R}, \bb{R}\times A
\in \mathfrak{S}_{n+1}$.
\item $\{(x_{1},...,x_{n})\in \bb{R}^{n}| x_{i}=x_{j}\} \in 
\mathfrak{S}_{n}$ for $1\leq i<j\leq n$.
\item $A\in \mathfrak{S}_{n+1} \Rightarrow \pi(A)\in \mathfrak{S}_{n}$, where
$\pi: \bb{R}^{n+1}\to \bb{R}^{n}$ is the projection on the first coordinates.
\item $\{r\}\in \mathfrak{S}_{1}$ for each $r\in \bb{R}$, and $\{(x,y)\in
\bb{R}^{2}|x<y\}\in \mathfrak{S}_{2}$.
\item The only sets in $\mathfrak{S}_{1}$ are the finite unions of open 
intervals and points.
\end{enumerate}

This notion is an elegant generalization of semialgebraic geometry (a 
standard reference for this theory is \cite{BCR}) and the semilinear 
or semialgebraic sets give the simplest examples of an o-minimal structure.\\

$\mathfrak{S}$ is called an o-minimal structure on the real field
$(\bb{R},+,\cdot)$, if it contains
the graph of addition and multiplication on $\bb{R}$ (and therefore
all semialgebraic sets). \\

For the construction of other o-minimal structures, see \cite{KarMac,LiRo,PeSpSt,
Speiss,RoSpWi,vDrSp,vDrSp2,Wilkie} and the references in \cite{vDrMi}. 
We use the notation
$\bb{R}_{an}$ for the o-minimal structure of globally subanalytic sets 
(i.e. subsets
of $\bb{R}^{n}$, that are subanalytic as subsets of the larger projctive
space $\bb{P}^{n}(\bb{R})$).\\  

An abstraction of the theory of o-minimal structures on the real field 
to the more general context of real analytic manifolds is the theory of
'analytic geometric categories', as introduced and studied by van den Dries 
and Miller \cite{vDrMi}. This is a generalization of the theory of subanalytic 
subsets of real analytic manifolds.

An analytic geometric category $\mathcal{S}$ is given if each real analytic 
manifold $M$ is equipped with a collection $\mathcal{S}(M)$ of subsets of $M$
such that the following conditions are satisfied (for each such manifold):
\begin{enumerate}
\item[AG1.] $\mathcal{S}(M)$ is a boolean algebra of subsets of $M$, with
$M\in \mathcal{S}(M)$.
\item[AG2.] If $A\in \mathcal{S}(M)$, then $A\times\bb{R}\in \mathcal{S}
(M\times\bb{R})$.
\item[AG3.] If $f:M\to N$ is a proper analytic map and $A\in \mathcal{S}(M)$, 
then $f(A)\in \mathcal{S}(N)$.
\item[AG4.] If $A\subseteq M$ and $(U_{i})_{i\in I}$ is an open covering of 
$M$, then $A\in \mathcal{S}(M)$ if and only if $A\cap U_{i}\in \mathcal{S}
(U_{i})$ for all $i\in I$.
\item[AG5.] Every bounded set in $\mathcal{S}(\bb{R})$ has a finite boundary.
\end{enumerate}

$\mathcal{S}$ corresponds uniquely (\cite[D.10]{vDrMi}) to an o-minimal
structure $\mathfrak{S}$ on $\bb{R}_{an}$ (i.e. an o-minimal structure
on the real field containing the structure $\bb{R}_{an}$ of globally 
subanalytic subsets). $\mathfrak{S}$ is defined by the subsets of 
$\bb{R}^{n}$, which belong to $\mathcal{S}$ as a subset of $\bb{P}^{n}
(\bb{R})$ (for the standard inclusion $\bb{R}^{n} \hookrightarrow \bb{P}^{n}
(\bb{R})$). Moreover, $\mathcal{S}$ can be recovered as the subsets of real 
analytic manifolds, which are locally (at each point of the ambient manifold) 
real analytic isomorphic to  
$\mathfrak{S}$-sets. Note that the last step is well defined, since  
the real analytic isomorphisms between bounded open subanalytic
subsets of affine spaces are definable in $\mathfrak{S}$.\\

We recall in this paper our variant of this notion, by restricting 
ourselves to 
(analytic) Nash manifolds (so in the above constructions, one should only look
at Nash manifolds and semialgebraic maps). In this way, every
o-minimal structure $\mathfrak{S}$ on the real field $(\bb{R},+,\cdot)$ (and 
not only those on $\bb{R}_{an}$) corresponds uniquely to a category, which we
call a Nash geometric category. This notion is a generalization of the 
theory of locally semialgebraic subsets of Nash manifolds. \\

A Nash geometric category $\mathcal{S}$ is given if each (real analytic) 
Nash manifold $M$ is equipped with a collection $\mathcal{S}(M)$ of subsets 
of $M$ such that the conditions NG$i.:=$AG$i., i=1,2,5$ and the
following conditions are satisfied (for each such manifold):
\begin{enumerate}
\item[NG3.] If $f:M\to N$ is a proper Nash map (i.e. analytic with a
semialgebraic graph) and $A\in \mathcal{S}(M)$, 
then $f(A)\in \mathcal{S}(N)$.
\item[NG4.] If $A\subseteq M$ and $(U_{i})_{i\in I}$ is a covering of 
$M$ by open semialgebraic subsets, then $A\in \mathcal{S}(M)$ if and only if 
$A\cap U_{i}\in \mathcal{S}(U_{i})$ for all $i\in I$.
\end{enumerate}

Then all results of \cite{vDrMi} extend to this context, if one makes
the obvious modifications (which way state and explain in the next section). 
Especially, an 'analytic geometric category' induces by 
restriction a 'Nash geometric category', and this correspondence is injective,
since the associated o-minimal structure is unique. A Nash geometric 
category $\mathcal{S}$ induces an o-minimal structure $\mathfrak{S}
(\mathcal{S})$ on 
$(\bb{R},+,\cdot)$, by the same definition as for analytic geometric
categories ($\bb{R}^{n}$ is an open semialgebraic subset of the Nash
manifold $\bb{P}^{n}(\bb{R})$). Moreover, $\mathcal{S}$ can be recovered as 
the subsets of Nash manifolds, which are locally (at each point of the 
ambient manifold) Nash isomorphic to  
$\mathfrak{S}$-sets. Note that the last step is well defined, since  
the Nash isomorphisms between bounded open semialgebraic subsets of affine 
spaces are definable in $\mathfrak{S}$. \\

The advantage of the notion
of a Nash geometric category is the fact, that this last step
is also well defined, if one starts with a category of 
$\mathfrak{X}$-sets, as defined by Shiota \cite{Shiota}, who uses a different
axiomatic setting.\\

A family of subsets of all affine spaces $\bb{R}^{n}$ is called 
$\mathfrak{X}$ (\cite[p.viii, p.95,96]{Shiota}), if it satisfies the axioms:
\begin{enumerate}
\item[$\mathfrak{X}$(i)] Every algebraic set in any Euclidean space is an 
element of $\mathfrak{X}$.
\item[$\mathfrak{X}$(ii)] If $X_{1}\subset \bb{R}^{n}$ and $X_{2}\subset 
\bb{R}^{n}$ are elements of $\mathfrak{X}$, then $X_{1}\cap X_{2},X_{1}
\backslash X_{2}$ and $X_{1}\times X_{2}$ are elements of $\mathfrak{X}$.
\item[$\mathfrak{X}$(iii)] If $X\subset \bb{R}^{n}$ is an element of 
$\mathfrak{X}$ and $p: \bb{R}^{n}\to \bb{R}^{m}$ is a linear map such that the 
restriction of $p$ to $\bar{X}$ is proper, then $p(X)$ is an element of 
$\mathfrak{X}$.
\item[$\mathfrak{X}$(iv)] If $X\subset \bb{R}$ and $X\in \mathfrak{X}$, then
each point of $X$ has a neigborhood in $X$, which is a finite 
union of points and intervals.
\end{enumerate}

Note that this a generalization of the notion of an o-minimal structure
on $(\bb{R},+,\cdot)$ (which is the same as a family $\mathfrak{X}_{0}$ in
the notation of Shiota). Moreover, a Nash isomorphisms $h$ between bounded 
open semialgebraic subsets $U_{i}$ ($i=1,2$) of affine spaces is definable in
$\mathfrak{X}$, and $X\subset U_{1}$ belongs to $\mathfrak{X}$ if and only if
$h(X)\subset U_{2}$ belongs to $\mathfrak{X}$ (by axiom $\mathfrak{X}$(iii)
and the graph embedding, as in \cite[II1.6]{Shiota}). 
So we get a well defined notion of subsets of (analytic) Nash manifolds,
which are locally Nash isomorphic to some $\mathfrak{X}$-set.
Note, that this argument doesn't work, if we do not use 'bounded charts'
(so it doesn't make sense, to introduce the notion of $\mathfrak{X}$-sets
of Nash manifolds. This notion works only for an o-minimal structure on the 
real field, compare \cite[p.507/508]{vDrMi}).\\
  
Sometimes, Shiota assumes in addition the axiom (\cite[p.97]{Shiota}):
\begin{enumerate}
\item[$\mathfrak{X}$(v)] If a subset $X$ of $\bb{R}^{n}$ is an 
$\mathfrak{X}$-set locally at each point of $\bb{R}^{n}$, then $X$ is an
$\mathfrak{X}$-set.
\end{enumerate}

Note that an o-minimal structure never satisfies this axiom. But if we 
restrict a Nash geometric category $\mathfrak{S}$ to all affine spaces 
$\bb{R}^{n}$, we get a family $\mathfrak{X}$ satisfying this axiom. 
Indeed, the notion of a Nash geometric category is equivalent to a family 
$\mathfrak{X}$, satisfying the axiom $\mathfrak{X}$(v). More precisely, the 
arguments of \cite{vDrMi} imply the following (see the next section): 

\begin{thm} \label{thm:comp}
\begin{enumerate}
\item A family $\mathfrak{X}$ induces a Nash geometric category $\mathcal{S}
(\mathfrak{X})$, whose sets are the subsets of Nash manifolds, which are 
locally (at each point of the ambient manifold) Nash-isomorphic to 
$\mathfrak{X}$-sets.
\item A Nash geometric category $\mathcal{S}$ induces: 
\begin{enumerate} 
\item[(a)] An o-minimal structure $\mathfrak{S}(\mathcal{S})$ on $(\bb{R},+,
\cdot)$, whose sets are
the subsets of $\bb{R}^{n}$, which belong as a subset of $\bb{P}^{n}(\bb{R})$
(for the standard inclusion $\bb{R}^{n} \hookrightarrow \bb{P}^{n}(\bb{R})$) 
to $\mathcal{S}$.
\item[(b)] By restriction to $\bb{R}^{n}$ ($n\in \bb{N}$) a family 
$\mathfrak{X}(\mathcal{S})$ satisfying the axiom $\mathfrak{X}(v)$.
\end{enumerate} 
\item $\mathcal{S}(\mathfrak{S}(\mathcal{S})) = \mathcal{S} = 
\mathcal{S}(\mathfrak{X}(\mathcal{S}))$.\\
\item $\mathfrak{S}(\mathcal{S}(\mathfrak{X})) \subseteq 
\mathfrak{X} \subseteq \mathfrak{X}(\mathcal{S}(\mathfrak{X}))$,
with $\mathfrak{S}(\mathcal{S}(\mathfrak{X})) = \mathfrak{X}$, if 
$\mathfrak{X}$ is an o-minimal structure on $(\bb{R},+,\cdot)$,
and $\mathfrak{X} = \mathfrak{X}(\mathcal{S}(\mathfrak{X}))$, if
$\mathfrak{X}$ satisfies the axiom $\mathfrak{X}(v)$.
\end{enumerate} 
\end{thm}

Remark, that this theorem allows an easy comparison between the different 
notions of 'geometric categories':\\

We get a one to one correspondence between o-minimal structures 
$\mathfrak{S}$ on $(\bb{R},+,\cdot)$, Nash geometric categories $\mathcal{S}$
and families $\mathfrak{X}$, satisfying the axiom $\mathfrak{X}$(v).   
Moreover, the analytic geometric categories (as a subset of the Nash geometric 
categories) correspond in this way to the o-minimal structures on 
$\bb{R}_{an}$, and to the families $\mathfrak{X}$, containing the subanalytic 
subsets of affine spaces and satisfying the axiom $\mathfrak{X}$(v)
(since the above correspondence is compatible with the natural partial order
on the families, which is induced by the inclusions of the subsets of the 
families).\\

Notice that $\mathfrak{S}(\mathcal{S}(\mathfrak{X}))$ is the greatest 
o-minimal structure containd in $\mathfrak{X}$. Every bounded
(especially every compact) $\mathfrak{X}$-set belongs to it, and every
$\mathfrak{X}$-map between such sets is definable in this structure. 
More precisely, the proof of theorem \ref{thm:comp} gives the following
characterization:
\begin{enumerate} \label{char:*}
\item[(*)] For each $A\subseteq \bb{R}^{n}$, $A\in \mathfrak{S}(\mathcal{S}
(\mathfrak{X}))$ if and only if $\tau_{n}(A)$ is an $\mathfrak{X}$-set
locally at each point of $\bb{R}^{n}$.
\end{enumerate} 

Here we use the semialgebraic map $\tau_{n}: \bb{R}^{n}\to 
\bb{R}^{n}$, given by 
\[\tau_{n}(x_{1},...,x_{n}):=(x_{1}(1+x_{1}^{2})^{-1/2}),...,
x_{n}(1+x_{n}^{2})^{-1/2}).\]
Note that $\tau_{n}$ is a Nash isomorphism of $\bb{R}^{n}$ onto 
$]-1,1[^{n}$.\\   

Similarly $\mathfrak{X}(\mathcal{S}(\mathfrak{X}))$ is the smallest family of 
this type, which contains $\mathfrak{X}$ and satisfies the axiom 
$\mathfrak{X}$(v). Moreover, it can be directly constructed as the family of 
sets, which are locally an $\mathfrak{X}$-set at each point of the ambient 
affine space (compare \cite[p.268]{Shiota}). 
For most of the results of \cite{Shiota} one can therefore assume that one
works with an o-minimal structure, or with a family $\mathfrak{X}$ 
satisfying the axiom $\mathfrak{X}$(v).
Moreover, one can extend these results into the framework of analytic
or Nash geometric categories, by using a suitable affine embedding of the 
ambient manifolds.  

\section{Geometric categories} 
In this section, we prove theorem \ref{thm:comp} and explain the modifications
compared to \cite{vDrMi}, which are necessary to translate their results into
the context of Nash geometric categories. We assume that the reader is familar
with \cite{vDrMi}. We follow the notation and numbering of 
\cite[sec.1,app.D]{vDrMi}.\\

For example, the fact that the family of subanalytic subsets of a real analytic
manifold is the 'smallest' analytic geometric category (\cite[p.499]{vDrMi}), 
translates of course into the fact that the family of locally semialgebraic 
subsets of Nash manifolds is the 'smallest' Nash geometric category.

Fix a Nash geometric category $\mathcal{S}$, let $M,N$ be (analytic) Nash
manifolds, and let $A\in \mathcal{S}(M)$,$B\in \mathcal{S}(N)$.
\begin{enumerate}
\item[1.1'] Every Nash map $f: M\to N$ (i.e. f is 
analytic with semialgebraic graph) is an $\mathcal{S}$-map.
\item[1.2'] Given an covering $(U_{i})$ of $M$ by open semialgebraic 
subsets, a map $f: A\to N$ is an $\mathcal{S}$-map if and only if each 
restriction
$f|U_{i}\cap A: U_{i}\cap A\to N$ is an $\mathcal{S}$-map.
\end{enumerate}
Moreover, the statements $1.i$ for $i=3,..,20$ of \cite{vDrMi} remain true
in the Nash geometric context without any modification (for example $1.7$ 
implies $cl(A):=\bar A\in \mathcal{S}(M)$).

For their proof, we modify \cite[app.D]{vDrMi} in the following way.
\begin{enumerate}
\item[D.1'] All locally semialgebraic sets are $\mathcal{S}$-sets;
in particular $\bb{R}^{n}\in \mathcal{S}(\bb{P}^{n}(\bb{R}))$. 
\end{enumerate}

Since this is a local statement, it suffices to show that the sets $\{x\in M|
f(x)=0\}$ and $\{x\in M|f(x)>0\}$ belong to $\mathcal{S}$(M), if $M$ is an 
open semialgebraic subset of some affine space and $f: M\to \bb{R}$ is a 
polynomial map. Then the proof of \cite[p.530]{vDrMi} applies (because the
embeddings $M\to M\times \bb{P}^{1}(\bb{R}), x\mapsto (x,f(x)),
x\mapsto (x,0)$ and the projection $M\times \bb{P}^{1}(\bb{R})\to M$ are 
proper Nash maps, and $\{(x,y)\in M\times \bb{P}^{1}(\bb{R})| f(x)\cdot 
y_{1}^{2}-y_{2}^{2}=0\}$ is locally given by the vanishing of a 
polynomial).\\

This implies also D.2'$:=$1.1', since the graph of a Nash map is (locally) 
semialgebraic, and D.3'$:=$1.2' follows directly from the axiom NG4.
Moreover, the statements $D.i$ for $i=4,..,9$ of \cite{vDrMi} remain true
in the Nash geometric context, with the same proof. For later applications let
us just recall that $D.4$ implies the stability of $\mathcal{S}$-sets under
products, und $D.6$ implies the stability under proper images:\\
Let $A,A'\in \mathcal{S}(M)$ with $A'\subseteq A$, $A$ locally closed and let
$f: A\to N$ be a proper $\mathcal{S}$-map. Then $f(A')\in \mathcal{S}(N)$.\\

Now we come to the proof of theorem \ref{thm:comp} (which is a generalization
of \cite[D.10]{vDrMi}).
\begin{proof}
(1) The statement 1. corresponds to \cite[D.10.3]{vDrMi}.
We already explained in the introduction, that $\mathcal{S}(\mathfrak{X})$
is well defined by the axioms $\mathfrak{X}$(i-iii). Moreover, NG4. follows 
from the definiton and NG.1,2 follow from $\mathfrak{X}$(i,ii).
For NG3., we start with a proper Nash map $f:M\to N$ and an $\mathcal{S}
(\mathfrak{X})$-set $A$ in $M$. Then the proof of \cite[p.533,534]{vDrMi}
applies, if we take for $y\in N$ an open semialgebraic neigborhood $V$ of $y$,
with a Nash isomorphism $h: V\to h(V)$ onto an open bounded (!) semialgebraic
subset of $\bb{R}^{n}$ containing $[-1,1]^{n}$ (and similarly for $g_{x}: U_{x}
\to g_{x}(U_{x})$, with $x\in f^{-1}h^{-1}([-1,1]^{n}), f(U_{x})\subset V$ and
such that $g_{x}(A\cap U_{x})$ belongs to $\mathfrak{X}$).
Then the map 
\[g_{x}(U_{x})\to h(V),\; a\mapsto h(f(g_{x}^{-1}(a)))\] 
is a Nash map between open bounded (!) semialgebraic subsets of affine spaces.
\[g_{x}(A\cap U_{x})\cap[-1,1]^{n} \subset g_{x}(U_{x})\] 
belongs to 
$\mathfrak{X}$ by $\mathfrak{X}$(i,ii). Therefore its image under the above
map belongs to $\mathfrak{X}$ by $\mathfrak{X}$(iii) (and the graph 
embedding). The rest of the proof works without changes and uses then only  
$\mathfrak{X}$(i,ii). Similarly, the proof of AG5. in \cite[p.534]{vDrMi}
applies also to the proof of NG5. in our situation, and uses 
$\mathfrak{X}$(iv) (instead of condition (6) for o-minimal structures).\\

(2) The statement 2.(a) corresponds to \cite[D.10.2]{vDrMi} and their
proof applies without changes. The proof of 2.(b) goes as follows:\\
$\mathfrak{X}$(i) follows from D.1' and $\mathfrak{X}$(ii) follows from
NG1. and D.4 ('stability under products'). Moreover, every bounded
$\mathcal{S}$-set in $\bb{R}$ belongs by 2.(a) to an o-minimal structure
and is therefore a finite union of intervals and points. This implies
$\mathfrak{X}$(iv) (by D.1' and NG1.). With $X$ belongs also the closure
$\bar{X}$ to $\mathcal{S}$ (this is a special case of 1.7, which will follow
from the first part of the statement 3. of the theorem). Then D.2'$:=$1.1' and 
D.6 imply the condition $\mathfrak{X}$(iii).\\

(3) The first equality in the statement 3. corresponds to the first part of
\cite[D.10.3]{vDrMi}. The proof in \cite[p.534]{vDrMi} applies (if one works 
with Nash isomorphisms instead of analytic isomorphisms).\\
Let $A\in \mathcal{S}(\mathfrak{X}(\mathcal{S}))(M)$. By definition,
there exists for all $x\in M$ a Nash isomophism $h_{x}:U_{x}\to V_{x}$, with 
$U_{x}$ an open semialgebraic neighborhood of $x$ 
in $M$ and $V_{x}$ an open semialgebraic subset      
of some affine space such that $h_{x}(A\cap U_{x})$ belongs to $\mathcal{S}
(V_{x})$. But then $A\cap U_{x}$ belongs to $\mathcal{S}(U_{x})$ (by NG3.)
for all $x\in M$, and therefore $A\in \mathcal{S}(M)$ (by NG4.).

Altogether we get $\mathcal{S}(M)\subseteq \mathcal{S}(\mathfrak{S}
(\mathcal{S}))(M)\subseteq \mathcal{S}(\mathfrak{X}(\mathcal{S}))(M)
\subseteq \mathcal{S}(M)$.\\

(4) First note, that \cite[D.10.1]{vDrMi} gives the  
characterization (*) of $\mathfrak{S}(\mathcal{S}(\mathfrak{X}))$-sets (the 
proof applies without
changes, since $\tau_{n}$ has a semialgebraic graph and $N\times \bb{P}^{n}
(\bb{R})\to N$ is a proper Nash map for each Nash manifold $N$).
Here we also use the obvious fact, that $A\subset \bb{R}^{n}$ belongs to
$\mathcal{S}(\mathfrak{X})$, if and only if $A$ is an $\mathfrak{X}$-set
locally at each point of $\bb{R}^{n}$.
This fact implies already the second part of statement 4.

Let $A\subset \bb{R}^{n}$ belong to $\mathfrak{X}$. Then $\tau_{n}(A)$ 
belongs to $\mathfrak{X}$ and therefore, $A$ belongs to 
$\mathfrak{S}(\mathcal{S}(\mathfrak{X}))$ by condition (*), if
$\mathfrak{X}$ is an o-minimal structure (this argument breaks down for 
more general $\mathfrak{X}$). Conversely, suppose $A\subset \bb{R}^{n}$ 
belongs to $\mathfrak{S}(\mathcal{S}(\mathfrak{X}))$ so that $\tau_{n}(A)$
is an $\mathfrak{X}$-set locally at each point of $\bb{R}^{n}$.
For each $x\in \bb{R}^{n}$, there exists $U_{x}$ open in
$\bb{R}^{n}$ such that $\tau_{n}(A)\cap U_{x}$ belongs to 
$\mathfrak{X}$. Since $\tau_{n}(A)$ is bounded, finitely many of the 
$U_{x}$ cover
$\tau_{n}(A)$ so that $\tau_{n}(A)$ belongs to $\mathfrak{X}$  
(by $\mathfrak{X}$(ii)). Hence, $A$ belongs to $\mathfrak{X}$ (by 
$\mathfrak{X}$(iii), since $\tau_{n}(A)$ is bounded. Compare with 
\cite[II.1.6]{Shiota}). This proves the first part of statement 4., and the
proof of theorem \ref{thm:comp} is finished.
\end{proof}

We finally remark, that the statements
$D.i$ for $i=11,..,17,19$ of \cite{vDrMi} remain true
in the Nash geometric context, with the same proof
(if one uses in the proof of D.11 (or D.19) an open covering $(U_{i})_{i\in
N}$ with $U_{i}$ open semialgebraic subsets of $M$ and $\phi_{i}:
U_{i}\to \bb{R}^{m}$ (or $h_{i}: U_{i}\to \bb{R}^{m}$) analytic Nash 
isomorphisms (with $h_{i}(U_{i}\cap A)\in \mathfrak{S}_{m}$).

\end{document}